\documentclass[12pt]{amsart}

\usepackage{amsthm}
\usepackage{amsmath}
\usepackage{amssymb}
\usepackage{eucal}
\usepackage{amsfonts}
\usepackage{amsmath}
\usepackage[Lenny]{fncychap}
\usepackage{enumerate}
\usepackage{amssymb}
\usepackage[english]{babel}
\usepackage[latin1]{inputenc}
\usepackage[breaklinks=true]{hyperref}
\usepackage[all]{xy}
\usepackage{fancybox}
\usepackage{latexsym,mathrsfs,longtable,lscape,dsfont,yfonts}
\usepackage{srcltx}
\usepackage{soul}
\usepackage{tikz}

\ifx\pdfpagewidth\undefined 
\usepackage[T1]{fontenc}
\else 
\usepackage[OT1]{fontenc} 
\fi 
\usepackage{amsfonts,amssymb,latexsym,longtable,amsmath}

\newtheorem{thm}{Theorem}[section]
\newcommand{\bthm}{\begin{thm}} \newcommand{\ethm}{\end{thm}}
\newtheorem{prop}[thm]{Proposition}
\newcommand{\bprp}{\begin{prop}} \newcommand{\eprp}{\end{prop}}
\newtheorem{fact}[thm]{Fact}
\newcommand{\bfct}{\begin{fact}} \newcommand{\efct}{\end{fact}}
\newtheorem{prob}[thm]{Problem}
\newcommand{\bprb}{\begin{prob}} \newcommand{\eprb}{\end{prob}}
\newtheorem{quest}[thm]{Question}
\newcommand{\bqtn}{\begin{quest}} \newcommand{\eqtn}{\end{quest}}
\newtheorem{lem}[thm]{Lemma}
\newcommand{\blem}{\begin{lem}} \newcommand{\elem}{\end{lem}}
\newtheorem{claim}[thm]{Claim}
\newcommand{\bclm}{\begin{claim}} \newcommand{\eclm}{\end{claim}}
\newtheorem{cor}[thm]{Corollary}
\newcommand{\bcor}{\begin{cor}} \newcommand{\ecor}{\end{cor}}
\newtheorem{conj}[thm]{Conjecture}
\newcommand{\bcnj}{\begin{conj}} \newcommand{\ecnj}{\end{conj}}
\theoremstyle{definition}
\newtheorem{defn}[thm]{Definition}
\newcommand{\bdfn}{\begin{defn}} \newcommand{\edfn}{\end{defn}}
\newtheorem{spec}[thm]{Specializing}
\newcommand{\bspc}{\begin{spec}} \newcommand{\espc}{\end{spec}}
\newtheorem{rem}[thm]{Remark}
\newcommand{\brem}{\begin{rem}} \newcommand{\erem}{\end{rem}}
\newtheorem{cnv}[thm]{Convention}
\newcommand{\bcnv}{\begin{cnv}} \newcommand{\ecnv}{\end{cnv}}
\theoremstyle{remark}
\newtheorem{exam}[thm]{Example}
\newcommand{\bexm}{\begin{exam}} \newcommand{\eexm}{\end{exam}}
\newcommand{\bpf}{\begin{proof}} \newcommand{\epf}{\end{proof}}

\newtheorem*{rep@theorem}{\rep@title}
\newcommand{\newreptheorem}[2]{%
\newenvironment{rep#1}[1]{%
 \def\rep@title{#2 \ref{##1}}%
 \begin{rep@theorem}}%
 {\end{rep@theorem}}
 }
\makeatother

\newtheorem{thmy}{\textbf{Theorem}}
\newenvironment{thmx}{\stepcounter{thm}\begin{thmy}}{\end{thmy}}

\newreptheorem{theorem}{Theorem}

\newreptheorem{corollary}{Corollary}

\newcommand{\R}{\mathbb R}

\newcommand{\lra}{\longrightarrow}

\renewcommand{\phi}{\varphi}
\renewcommand{\theta}{\vartheta}

\newcommand{\gep}{{\epsilon}}

\newcommand{\gd}{{\delta}}

\newcommand{\gb}{{\beta}}
\newcommand{\gom}{{\omega}}

\newcommand{\mkp}{\medskip}
\newcommand{\bkp}{\bigskip}
\def\defi{\buildrel\rm def \over=}

\begin{document}

\title[Equicontinuity criteria for metric-valued sets of continuous functions]{Equicontinuity criteria for metric-valued sets of continuous functions}

\author[M. Ferrer]{Marita Ferrer}
\address{Universitat Jaume I, Instituto de Matem\'aticas de Castell\'on,
Campus de Riu Sec, 12071 Castell\'{o}n, Spain.}
\email{mferrer@mat.uji.es}

\author[S. Hern\'andez]{Salvador Hern\'andez}
\address{Universitat Jaume I, INIT and Departamento de Matem\'{a}ticas,
Campus de Riu Sec, 12071 Castell\'{o}n, Spain.}
\email{hernande@mat.uji.es}

\author[L. T\'arrega]{Luis T\'arrega}
\address{Universitat Jaume I, IMAC and Departamento de Matem\'{a}ticas,
Campus de Riu Sec, 12071 Castell\'{o}n, Spain.}
\email{ltarrega@uji.es}

\thanks{ Research Partially supported Universitat Jaume I, grant P1-1B2015-77.
The second
author acknowledges partial support by Generalitat Valenciana,
grant code: PROMETEO/2014/062, and the third author also acknowledges partial support of the Spanish Ministerio de Econom\'{i}a y Competitividad
grant MTM 2013-42486-P}

\begin{abstract}
Combining ideas of Troallic \cite{Troallic1996} and Cascales, Namioka, and Vera \cite{Cascales2000},
we prove several characterizations of \textit{almost equicontinuity} and
\textit{hereditarily almost equicontinuity} for subsets of metric-valued continuous functions
when they are defined on a \v{C}ech-complete space. We also obtain some applications
of these results to topological groups and dynamical systems.
\end{abstract}



\date{\today}
\maketitle \setlength{\baselineskip}{24pt}

\section{Introduction}
Let $X$ and $(M,d)$ be a Hausdorff, completely regular space and a metric space, respectively,
and let $C(X,M)$ denote the set of all continuous functions from $X$ to $M$.
A subset  $G\subseteq C(X,M)$ is said to be \textit{almost equicontinuous} if $G$ is equicontinuous on a dense subset of $X$.
If $G$ is almost equicontinuous for every closed nonempty subset of $X$, then it is said that $G$ is \textit{hereditarily almost equicontinuous}.
{The main goal of this paper is to extend to arbitrary topological spaces these two important notions, which
were introduced in the setting of topological dynamics studying the enveloping semigroup of a flow \cite{Akin1998,Glasner2000,Glasner2006}.}

In addition to their intrinsic academic interest, it turns out that these two concepts have found application in other different settings as it will be
made clear in the sequel. First, we shall provide some basic notions and terminology.

{Given } $F\subseteq X$, the symbol $t_p(F)$ (resp. $t_{\infty}(F)$) will denote the topology, on $C(X,M)$,
of pointwise convergence (resp. uniform convergence) on $F$.
{For a set $G$ of functions from $X$ to $M$} and $Z\subseteq X$, the symbol $G\vert_Z$ will denote
the set $\lbrace g\vert_Z:g\in G\rbrace$. {We denote by $\overline{G}^{M^X}$ the closure of $G$ in the Tychonoff product space $M^X$.}\\
{The symbolism $(F,t_p(\overline{G}^{M^X}))$ will denote the set $F$
equipped with the weak topology generated by the functions in $\overline{G}^{M^X}\vert_F$.}
In like manner, the symbol $[A]^{\leq \omega}$ will denote the set of all countable subsets of $A$.
A topological space $X$ is said to be \emph{\v{C}ech-complete} if it is a $G_\delta$-subset of its Stone-\v{C}ech compatification.
The family of \v{C}ech-complete spaces  includes all complete metric spaces and all locally compact spaces.
Several quotient spaces are used along the paper. For the reader's sake, a detailed description of them is presented at the Appendix.

We now formulate our main results.

\begin{thmx}\label{resultado_simpli}
Let $X$ and $(M,d)$ be a \v{C}ech-complete space and a separable metric space, respectively, and let $G\subseteq C(X,M)$
such that $\overline{G}^{M^X}$ is compact.
Consider the following three properties:
\begin{enumerate}[(a)]
\item $G$ is almost equicontinuous.
\item There exists a dense Baire subset $F\subseteq X$ such that $(\overline{G}^{M^X})\vert_F$ is metrizable.
\item There exists a dense $G_{\delta}$ subset $F\subseteq X$ such that
$(F,t_p(\overline{G}^{M^X}))$ is Lindel\"of.
\end{enumerate}
Then $(b)\Rightarrow (c) \Rightarrow (a)$.
If $X$ is also a hereditarily Lindel\"of space, then all conditions are equivalent.
\end{thmx}
\mkp


Next result characterizes hereditarily almost equicontinuous families of functions defined on a \v{C}ech-complete space  (this question has been
studied in detail in \cite{Ruiz2014} for compact spaces).

\begin{thmx}\label{resultado_6}
Let $X$ and $(M,d)$ be a \v{C}ech-complete space and a  metric space, respectively, and let $G\subseteq C(X,M)$
such that $\overline{G}^{M^X}$ is compact.
Then the following conditions are equivalent:
\begin{enumerate}[(a)]
\item $G$ is {hereditarily} almost equicontinuous.
\item $L$ is {hereditarily} almost equicontinuous on $F$, for all $L\in [G]^{\leq \omega}$ and $F$ a separable and compact subset of $X$.
\item $(\overline{L}^{M^X})\vert_F$ is metrizable, for all $L\in [G]^{\leq \omega}$ and $F$ a separable and compact subset of $X$.
\item $(F,t_p(\overline{L}^{M^X}))$ is Lindel\"{o}f, for all $L\in [G]^{\leq \omega}$ and $F$ a separable and compact subset of $X$.
\end{enumerate}
\end{thmx}

\brem If $G$ is a subset of $C(X,M)$ such that $K\defi \overline{G}^{M^X}$  is contained in $C(X,M)$, then the implication $(c)\Rightarrow (a)$ in
Theorem \ref{resultado_6} provides a different proof of the celebrated {Namioka Theorem \cite[Theorem~2.3]{Namioka1974a}.}
Indeed, given any $L\in [G]^{\leq \omega}$ and any  separable compact subset $F$ of $X$,
since $K\subseteq C(X,M)$ and $F$ is separable, it follows that $((\overline{L}^{M^X})\vert_F,t_p(F))$ is metrizable.
Thus $G$ (and therefore $K$) is hereditarily almost equicontinuous.
\erem

\bcor\label{nota_simp}
With the same hypothesis of Theorem \ref{resultado_6}, consider the following three properties:
\begin{enumerate}[(a)]
\item $G$ is hereditarily almost equicontinuous.
\item $G$ is {hereditarily} almost equicontinuous on $F$, for all $F$ a separable and compact subset of $X$.
\item $(F,t_p(\overline{G}^{M^X}))$ is Lindel\"{o}f, for all $F$ a separable and compact subset of $X$.
\end{enumerate}
Then $(a)\Leftrightarrow (b)\Leftarrow (c)$.
\ecor



\section{Applications}

The results formulated in the previous section have consequences in different settings.
First, we consider an application to fragmentability.

A topological space $X$ is said to be \textit{fragmented by a pseudometric $\rho$} if for each nonempty subset $A$ of $X$ and for each
 $\epsilon>0$ there exists a nonempty open subset $U$ of $X$ such that $U\cap A\neq \emptyset$ and $\rho$-$diam(U\cap A)\leq \epsilon$.
 This notion 
 was introduced by Jayne and Rogers in \cite{Jayne1985}. Further work has been done by many workers. It will suffice to mention here the contribution by
 Namioka \cite{Namioka1987} and Ribarska \cite{Ribarska1987}.

Let $X$ be a topological space, $(M,d)$ a metric space and $G\subseteq M^X$ a family of functions.
Whenever feasible, {for example if $\overline{G}^{M^X}$ is compact,} we will consider the pseudometric $\rho_{G,d}$, defined as follows:
\[ \rho_{G,d}(x,y)\defi\sup\limits_{g\in G}d(g(x),g(y)),\quad \forall x,y\in X.\]

Therefore, taking into account Definition \ref{def_AE} and Lemma \ref{lem1}, we have the following proposition.

\bprp\label{HAE_FRAG}
Let $X$ and $(M,d)$ be a topological space and a metric space, respectively, and let $G\subseteq C(X,M)$ such that $\overline{G}^{M^X}$ is compact. Consider the following two properties:
\begin{enumerate}[(a)]
\item $G$ is hereditarily almost equicontinuous.
\item $X$ is fragmented by $\rho_{G,d}$.
\end{enumerate}
Then $(a)$ implies $(b)$. If $X$ is a hereditarily Baire space, then $(a)$ and $(b)$ are equivalent.
\eprp

{As a consequence, we have the following corollary of Theorem B.}

\bcor\label{resultado_6_fragmentab}
Let $X$ and $(M,d)$ be a \v{C}ech-complete space and a metric space, respectively, and let $G\subseteq C(X,M)$
such that $\overline{G}^{M^X}$ is compact.
Then the following conditions are equivalent:
\begin{enumerate}[(a)]
\item $X$ is fragmented by $\rho_{G,d}$.
\item $F$ is fragmented by $\rho_{L,d}$, for all $L\in [G]^{\leq \omega}$ and $F$ a separable and compact subset of $X$.
\item $((\overline{L}^{M^X})\vert_F,t_p(F))$ is metrizable, for all $L\in [G]^{\leq \omega}$ and $F$ a separable and compact subset of $X$.
\item $(F,t_p(\overline{L}^{M^X}))$ is Lindel\"{o}f, for all $L\in [G]^{\leq \omega}$ and $F$ a separable and compact subset of~$X$.
\end{enumerate}
\ecor

It is easy to check that, in the context of {topological groups}, the notion of
almost equicontinuity is equivalent to equicontinuity.
This fact allows us to characterize equicontinuous subsets of group homomorphisms using Theorem A.

From here on, if $X$ and $M$ are topological groups, the symbol $CHom(X,M)$ will denote the set of continuous homomorphisms of $X$ into $M$. Recall that a topological group $G$ is said to be $\omega$-\textit{narrow} if for every neighborhood $V$ of the neutral element,
there exists a countable subset $E$ of $G$ such that $G=EV$.

\bcor\label{resultado_simpli_gt}
Let $X$ and $(M,d)$ be a \v{C}ech-complete topological group and a metric separable group, respectively, and let $G$ be a
subset of $CHom(X,M)$ such that  $\overline{G}^{M^X}$ is compact.
Consider the following three properties:
\begin{enumerate}[(a)]
\item $G$ is equicontinuous.
\item $G$ is relatively compact in $CHom(X,M)$ with respect to the compact open topology.
\item There exists a dense Baire subset $F\subseteq X$ such that $(\overline{G}^{M^X})\vert_F$ is metrizable.
\item There exists a dense $G_{\delta}$ subset $F\subseteq X$ such that
$(F,t_p(\overline{G}^{M^X}))$ is Lindel\"of.
\end{enumerate}
{Then }$(c)\Rightarrow (d) \Rightarrow (a)\Leftrightarrow (b)$.
If $X$ is also $\omega$-narrow, then all conditions are equivalent. Furthermore $(c)$ and $(d)$ are also true for $F=X$.
\ecor
\bpf
The equivalence $(a)\Leftrightarrow (b)$ follows from Ascoli Theorem. So, after Theorem A, it will suffice to show
the implication $(a)\Rightarrow (c)$ {for an $\omega$-narrow $X$}. Now, assuming that $G$ is equicontinuous, it follows that
$K\defi\overline{G}^{M^X}\subseteq CHom(X,M)$. Thus $K$ is an equicontinuous compact subset of continuous group homomorphisms.
As a consequence, it is known that $K$ is metrizable. (see \cite[Cor. 3.5]{Ferrer2012}).
\epf

Extending a result given by Troallic in \cite[Corollary~3.2]{Troallic1996}, we can reduce the ve\-ri\-fi\-cation
of hereditarily almost equicontinuity to countable subsets. The equivalence $(a)\Leftrightarrow (b)$
bellow is due to Troallic (op. cit.).

\bcor\label{resultado_6_cor}
Let $X$ and $(M,d)$ be a \v{C}ech-complete topological group and a metric group, respectively, and let $G$ be a
subset of $CHom(X,M)$ such that  $\overline{G}^{M^X}$ is compact.
Then the following conditions are equivalent:
\begin{enumerate}[(a)]
\item $G$ is equicontinuous.
\item $L$ is equicontinuous on $F$, for all $L\in [G]^{\leq \omega}$ and $F$ a separable and compact subset of $X$.
\item $((\overline{L}^{M^X})\vert_F,t_p(F))$ is metrizable, for all $L\in [G]^{\leq \omega}$ and $F$ a separable and compact subset of $X$.
\item $(F,t_p(\overline{L}^{M^X}))$ is Lindel\"{o}f, for all $L\in [G]^{\leq \omega}$ and $F$ a separable and compact subset of~$X$.
\end{enumerate}
\ecor
\mkp

{For a function $f:X\times Y\rightarrow M$ let $f_x:Y\rightarrow M$ ($f^y:X\rightarrow M$ ) be $f(x,\cdot)$
for a fixed $x\in X$ ($f(\cdot,y)$ for a fixed $y\in Y$, resp.).}

A variation of the celebrated Namioka Theorem \cite{Namioka1974a} is also obtained as a corollary of Theorems
\ref{resultado_simpli} and \ref{resultado_6} (cf. \cite{Kenderov2006,Talagrand1979,Piotrowski1985}).

\bcor\label{resultado_simpli_Namioka}
Let $X$, $H$, and $(M,d)$ be a \v{C}ech-complete space, a compact space, and a metric space, respectively, and
let $f:X\times H\rightarrow M$ be a map satisfying that $f_x\in C(H,M)$ for every $x\in X$ and there is
a dense subset $G$ of $H$ such that $f^g\in C(X,M)$ for every $g\in G$.
{Suppose that} any of the two following equivalent conditions holds.
\begin{enumerate}[(a)]
\item There exists a dense Baire subset $F\subseteq X$ such that $(\overline{G}^{M^X})\vert_F$ is metrizable.
\item There exists a dense $G_{\delta}$ subset $F\subseteq X$ such that
$(F,t_p(\overline{G}^{M^X}))$ is Lindel\"of.
\end{enumerate}
Then there exists a $G_{\delta}$ and dense subset $F$ in $X$ such that $f$ is jointly continuous at each point of $F\times H$.
\ecor

Finally, we obtain some applications to dynamical systems \cite{Glasner2000,Glasner2004,Glasner2006}. Recall that a dynamical system, or a $G$-space, is a Hausdorff space
$X$ on which a topological group $G$ acts continuously. We denote such a system by $(G,X)$. For each $g\in G$ we have the
self-homeomorphism $x\mapsto gx$ of $X$ that we call $g$-translation.

\bcor\label{resultado_simpli_sist_din}
Let $X$ be a Polish $G$-space such that $\overline{G}^{X^X}$ is compact.
The following properties are equivalent:
\begin{enumerate}[(a)]
\item $G$ is almost equicontinuous.
\item There exists a dense Baire subset $F\subseteq X$ such that $(\overline{G}^{X^X})\vert_F$ is metrizable.
\item There exists a dense $G_{\delta}$ subset $F\subseteq X$ such that
$(F,t_p(\overline{G}^{X^X}))$ is Lindel\"of.
\end{enumerate}
\ecor

\bcor\label{resultado_6_sist_din}
Let $X$ be a {completely} metrizable $G$-space such that $\overline{G}^{X^X}$ is compact.
Then the following conditions are equivalent:
\begin{enumerate}[(a)]
\item $G$ is hereditarily almost equicontinuous.
\item $L$ is {hereditarily} almost equicontinuous on $F$, for all $L\in [G]^{\leq \omega}$ and $F$ a compact subset of $X$.
\item $((\overline{L}^{M^X})\vert_F,t_p(F))$ is metrizable, for all $L\in [G]^{\leq \omega}$ and $F$ a compact subset of $X$.
\item $(F,t_p(\overline{L}^{M^X}))$ is Lindel\"{o}f, for all $L\in [G]^{\leq \omega}$ and $F$ a compact subset of $X$.
\end{enumerate}
\ecor

In \cite[Problem 28]{Arkhangelskii1990}, Arkhangel'skii raises the following question:
Let $X$ be a Lindel\"of space and let $K$ be a compact subset of $(C(X),t_p(X))$.
Is it true that the tightness of $K$ is countable?
As far as we know, this question
is still open in ZFC. Here we provide a partial answer to Arkhangel'skii's question.

\bcor
Let $X$ be a Lindel\"of space and let $K$ be a compact subspace of $(C(X),t_p(X))$. If there is a a dense subset
$G\subseteq K$ such that $(X,t_p(G))$ is \v{C}ech-complete and hereditarily Lindel\"of, then $K$ is metrizable.
\ecor
\bpf
The proof of this result is consequence of Theorem B. Indeed,
remark that, if $F$ is a subset of $X$ that is closed in the $t_p(G)$-topology, then
$F$ will be \v{C}ech-complete and hereditarily Lindel\"of as well.
Moreover, since $G\subseteq K$,
it follows that $F$ is also closed in the $t_p(K)$-topology and, as a consequence, Lindel\"of.
Applying Corollary \ref{nota_simp} to the (compact) space $K$, which is equipped with the $t_p(X)$-topology,
it follows that $G$ is hereditarily almost equicontinuous on $X$. Since $(X,t_p(G))$ is \v{C}ech-complete and
hereditarily Lindel\"of, Proposition \ref{GMU} yields the metrizability of $K=\overline{G}^{\R^X}$.
\epf

\section{Basic results}

Within the setting of dynamical systems, the following definitions appear in \cite{Akin1998}.

\bdfn\label{def_AE}
Let $X$ and $(M,d)$ be a topological space and a metric space respectively, and let $G\subseteq C(X,M)$.
According to \cite{Akin1998}, we say that a point $x\in X$ is an \emph{equicontinuity point} of $G$ when for every $\gep >0$
there is a neighborhood $U$ of $x$ such that $diam(g(U))<\gep$ for all $g\in G$.
We say that $G$ is \emph{almost equicontinuous} when the subset of equicontinuity points of $G$ is dense in $X$.
Furthermore, it is said that $G$ is \emph{hereditarily almost equicontinuous} if $G\vert_{A}$ is almost equicontinuous for every
nonempty closed subset $A$ of $X$.
\edfn

The proof of the following lemma is known. However it is very useful in order to obtain subsets of continuous functions that are not almost equicontinuous.
We include its proof here for completeness sake.

\blem\label{lem1}
Let $X$ and $(M,d)$ be a topological space and a metric space respectively, and let $G\subseteq C(X,M)$. Consider the following two properties:
\begin{enumerate}[(a)]
\item $G$ is almost equicontinuous.
\item For every nonempty open subset $U$ of $X$ and $\epsilon>0$, there exists a nonempty open subset $V\subseteq U$ such that $diam(g(V))<\epsilon$ for all $g\in G$.
\end{enumerate}

Then (a) implies (b). If X is a Baire space, then (a) and (b) are equivalent. Furthermore, in this case, the subset of equicontinuity
points of $G$ is a dense $G_\delta$-set in $X$.
\elem
\bpf
That (a) implies (b) is obvious. Assume that $X$ is a Baire space and (b) holds. Given $\epsilon>0$ arbitrary, we consider the open set
$O_{\epsilon}\defi\bigcup\lbrace U\subseteq X:\quad U$ is a nonempty open subset $\wedge$  $diam(g(U))<\epsilon\quad\forall g\in G \rbrace$. By (b), we have that
$O_{\epsilon}$ is nonempty and dense in $X$. Since $X$ is Baire, taking $W\defi\bigcap\limits_{n<\omega}O_{\frac{1}{n}}$,
we obtain a dense $G_\delta$ subset which is the subset of equicontinuity points of $G$.
\epf

\brem
As a consequence of assertion (b) in Lemma \ref{lem1}, it follows that, when $X$ is a Baire space, a subset of functions $G$ is hereditarily almost
equicontinuous if, and only if, $G\vert_{A}$ is almost equicontinuous for every nonempty (non necessarily closed) subset $A$ of $X$. Since we mostly work with
Baire spaces here, we will make use of this fact in some parts along the paper.
\erem

Note that the set of equicontinuity points of a subset of functions $G$ is a $G_{\delta}$-set. Next corollary is a straightforward consequence of Lemma \ref{lem1}.
\bcor\label{cor_cantor}
Let $X$ and $(M,d)$ be a topological space and a metric space respectively, and let $G\subseteq C(X,M)$. Suppose there is an open basis $\mathcal{V}$
in $X$ and $\gep >0$ such that for every $V\in \mathcal{V}$, there is $g_V\in G$ with $diam(g_V(V))\geq\gep$. Then $G$ is not almost equicontinuous.
\ecor


Let $2^{\omega}$ be the Cantor space and let $2^{(\omega)}$ denote the set of finite sequences of $0$'s and $1$'s.
For a $t\in 2^{(\omega)}$, we designate by $\vert t \vert$ the length of $t$. For $\sigma\in 2^{\omega}$ and $n>0$ we write $\sigma\vert n$ to denote {$(\sigma(0),\ldots,\sigma(n-1))\in 2^{(\omega)}$}. If $n=0$ then $\sigma\vert 0\defi \emptyset$.

Applying Corollary \ref{cor_cantor}, it is easy to obtain subsets of continuous functions that are not almost equicontinuous.

\bexm\label{ex_cantor1}
Let $X=2^\gom$ be the Cantor space and let $G=\left\lbrace \pi_n\right\rbrace_{n<\gom}$ be the set of all projections of $X$ onto $\{0,1\}$. Then
$G$ is not almost equicontinuous.
\eexm
\bpf
Let $U\neq\emptyset$ be an open subset in $X$. Then, for some index $n<\gom$ we have $\pi_n(U)=\{0,1\}$, which implies $diam(\pi_n(U))>1/2$.
Therefore $G$ is not almost equicontinuous by {Corollary \ref{cor_cantor}}.
\epf
\mkp

The precedent result can be generalized in order to obtain a more general example of non-almost equicontinuous set of functions.
It turns out that this example is universal in a sense that will become clear along the paper.

\bexm\label{ex_cantor2}
Let $X=2^\gom$ be the Cantor space and let $(M,d)$ be a metric space. Let $\{U_t : t\in 2^{(\gom)}\}$ be the canonical open
basis of $X$. If $G=\lbrace g_t\rbrace_{t\in 2^{(\gom)}}$ is a set of continuous functions {on $X$ into $M$} satisfying that
$diam(g_t(U_t))\geq\gep$ for some fixed $\gep >0$ and all $t\in 2^{(\gom)}$, then $G$ is not almost equicontinuous.
\eexm
\mkp

Next result gives a sufficient condition for the equicontinuity of a family of functions. It extends a well known result
by Corson and Glicksberg \cite{Corson1970}. However,
we remark that the subset $F$ found in the lemma below can become the empty set if $Z$ is a first category subset of $X$.

\blem\label{lem_AlmostEquicontinuity}
Let $X$ and $(M,d)$ be a topological space and a separable metric space, respectively. If $G\subseteq C(X,M)$
and  $(\overline{G}^{M^X})\vert_Z$ is {metrizable and compact} for some dense subset $Z$ of $X$,
then there is a residual subset $F$ in $Z$ such that $G$ is equicontinuous at every point in $F$.
In case $Z$ is of second category in $X$, it follows that $F$ will be necessarily nonempty.
\elem
\bpf
Set $H\defi\overline{G}^{M^X}$ and consider the map $eval:X\rightarrow C(H,M)$, $x\mapsto eval_x$; defined by $eval_x(f)\defi f(x)$ for all $x\in X$ and $f\in H$.

{For simplicity's sake, the symbols $C_{t_p(G)}(H\vert_Z,M)$ and $C_{\infty}(H\vert_Z,M)$ will denote the space $C(H\vert_Z,M)$
equipped with the pointwise convergence $t_p(G)$ and the uniform convergence topology, respectively.}

Now set $\Phi$ such that the following diagram commutes
\[
\xymatrix{ Z \ar@{>}[rr]^{eval} \ar[dr]^{\Phi} & & C_{t_p(G)}(H\vert_Z,M)\ar[dl]_{id} \\ & C_{\infty}(H\vert_Z,M) &}
\]

\mkp

Remark that the evaluation map, $eval$, is continuous because $G\subseteq C(X,M)$. Since $H\vert_Z$ is $t_p(Z)$-compact and metrizable and $Z$ is dense in $X$,
it follows that $C_\infty(H\vert_Z,M)$ is separable and metrizable (see \cite[Cor. 4.2.18]{Engelking1989}). Therefore, for every $n<\gom$, there is a sequence of closed balls
$\{\overline{B}(u_i^{(n)},1/n) : i< \gom\}$ that covers $C_\infty(H\vert_Z,M)$. Furthermore, since $G$ is dense in $H$, we have that
each $\overline{B}(u_i^{(n)},1/n)$ is also closed in $C_{t_p(G)}(H\vert_Z,M)$. As a consequence
$K_{(i,n)}\defi\Phi^{-1}(\overline{B}(u_i^{(n)},1/n))=eval^{-1}(\overline{B}(u_i^{(n)},1/n))$ is closed in $Z$ for all $i,n<\omega$, because $eval$ is continuous.

We have that $Z\subseteq \bigcup\limits_{i<\gom} K_{(i,n)}$ for every $n<\omega$, so $Z\subseteq \bigcap\limits_{n<\gom}\bigcup\limits_{i<\gom} K_{(i,n)}$.
Observe that  $\bigcup\limits_{n<\gom}\bigcup\limits_{i<\gom} (K_{(i,n)}\setminus int_Z(K_{(i,n)}))$ is a set of first category in $Z$.
As a consequence
$$F\defi Z\setminus \bigcup\limits_{n<\gom}\bigcup\limits_{i<\gom} (K_{(i,n)}\setminus int_Z(K_{(i,n)}))$$
is a residual set in $Z$.

We now verify that $G$ is equicontinuous at each point $z\in F$. Let $z\in F$ and $\gep >0$ arbitrary.
Take $n_0<\gom$ such that $2/n_0<\gep$. Since $z\in F\subseteq \bigcap\limits_{n<\gom}\bigcup\limits_{i<\gom} K_{(i,n)}\subseteq \bigcup\limits_{i<\gom} K_{(i,n_0)}$ there is $i_0<\gom$ such that $z\in K_{(i_0,n_0)}$. We claim that $z\in int_Z(K_{(i_0,n_0)})$. Indeed, if we assume that $z\not\in int_Z(K_{(i_0,n_0)})$, then $z\in K_{(i_0,n_0)}\setminus int_Z(K_{(i_0,n_0)})$. Therefore, $z\in \bigcup\limits_{n<\gom}\bigcup\limits_{i<\gom} (K_{(i,n)}\setminus int_Z(K_{(i,n)}))$ and $z\not\in F$, which is a contradiction.

Since $z\in int_Z(K_{(i_0,n_0)})$ there is a nonempty open set $A$ in $X$ such that $int_Z(K_{(i_0,n_0)})=A\cap Z$. Note that $A\cap Z$ is dense on $A$ because $Z$ is dense in $X$. So, $z\in A=\overline{A\cap Z}^{A}\subseteq \overline{A\cap Z}^{X}$.

Let $a,b\in A\cap Z$. Then $\Phi(a)=eval_a,\Phi(b)=eval_b\in \overline{B}(u_{i_0}^{(n_0)},1/n_0)$. Consequently, $d(g(a),g(b))\leq 2/n_{0}$ for every $g\in G$. So, given $x,y\in A\subseteq \overline{A\cap Z}^{X}$ we have that $d(g(x),g(y))\leq 2/n_{0}$ for every $g\in G$. Then $diam(g(A))\leq 2/n_{0}<\epsilon$ for all $g\in G$.
\epf

\brem\label{nota_K}
Let $X$ be a topological space, $(M,d)$ be a metric space and $G$ be a subset of $C(X,M)$ that we consider equipped with
the pointwise convergence topology $t_p(X)$ in the sequel, unless otherwise stated.

Set
$$ K\defi\lbrace \alpha:M\rightarrow [-1,1]:\vert \alpha(m_1)-\alpha(m_2)\vert\leq d(m_1,m_2),\quad\forall m_1,m_2\in M\rbrace.$$

\noindent It is readily seen that $K$ is a compact subspace of $[-1,1]^M$.

Consider the evaluation map $\varphi:X\times G\rightarrow M$ defined by $\varphi(x,g)\defi g(x)$ for all $(x,g)\in X\times G$,
which is clearly separately continuous. The map $\varphi$ has associated a separately continuous map $f:X\times (G\times K)\rightarrow [-1,1] $
defined by $f(x,(g,\alpha))\defi\alpha(g(x))$ for all $(x,(g,\alpha))\in X\times ( G\times K)$.

 Set $$\nu:\overline{G}^{M^X}\times K\rightarrow [-1,1]^{X}$$
 \noindent defined by
$$\nu(h,\alpha)\defi\alpha\circ h\ \hbox{for all}\ h\in \overline{G}^{M^X}\ \hbox{and}\ \alpha\in K.$$

We claim that $\nu$ is continuous. Indeed, let $\lbrace (h_{\delta},\alpha_{\delta})\rbrace_{\delta\in \Delta}\subseteq \overline{G}^{M^X}\times K$ be a net that converges to $(h,\alpha)\in \overline{G}^{M^X}\times K$. Given $\epsilon>0$ and $x\in X$, then there exists $\delta_0\in \Delta$ such that $d(h_{\delta}(x),h_0(x))<\epsilon/2$ and $\vert \alpha_{\delta}(h_0(x))-\alpha_0(h_0(x))\vert<\epsilon/2$ for all $\delta>\delta_0$. Therefore, we have that $\vert \nu(h_0,\alpha_0)(x)-\nu(h_{\delta},\alpha_{\delta})(x)\vert=\vert \alpha_0(h_0(x))-\alpha_{\delta}(h_{\delta}(x))\vert\leq\vert \alpha_0(h_0(x))-\alpha_{\delta}(h_0(x))\vert+\vert \alpha_{\delta}(h_0(x))-\alpha_{\delta}(h_{\delta}(x)) \vert\leq\vert \alpha_0(h_0(x))-\alpha_{\delta}(h_0(x))\vert+d( h_0(x),h_{\delta}(x))<\epsilon$ for all $\delta>\delta_0$.

\noindent Since $G\subseteq C(X,M)$, we have that $\nu(G\times K)\subseteq C(X,[-1,1])$.

\erem

For $m_0\in M$, define $\alpha_{m_0}\in [-1,1]^M$ by $\alpha_{m_0}(m)\defi d(m,{m_0})$ for all $m\in M$. It is easy to check that $\alpha_{m_0}\in K$.

\blem\label{lemma_homeo_sugerencia}
Let $X$ be a topological space, $(M,d)$ a metric space and $G$ a subset of $C(X,M)$. Let $K$ and $\nu$ be the space and the map defined in Remark \ref{nota_K}. Then, for every subset $F$ of $X$, the identity map $id:(F,t_p(\overline{G}^{M^X}))\rightarrow (F,t_p(\nu(\overline{G}^{M^X}\times K)))$ is a homeomorphism.
\elem
\begin{proof}
Let $\lbrace x_{\delta}\rbrace_{\delta\in \Delta}\subseteq F$ be a net that $t_p(\overline{G}^{M^X})$-converges to $x$.
Since $\alpha$ is continuous, for any $(h,\alpha)\in \overline{G}^{M^X}\times K$, we have $\lim\limits_{\delta\in\Delta}\nu(h,\alpha)(x_{\delta})=\lim\limits_{\delta\in\Delta}\alpha(h(x_{\delta}))=\alpha(h(x))=\nu(h,\alpha)(x)$. So, $id$ is continuous.
Conversely, let $\lbrace x_{\delta}\rbrace_{\delta\in \Delta}\subseteq F$ be a net that $t_p(\nu(\overline{G}^{M^X}\times K))$-converges to $x_0\in F$. Given $h\in \overline{G}^{M^X}$ arbitrary, take $\alpha_{h(x_0)}\in K$. So, fixed $\epsilon>0$, there is $\delta_0\in\Delta$ such that $\epsilon> \vert \nu(h,\alpha_{h(x_0)})(x_{\delta})-\nu(h,\alpha_{h(x_0)})(x_0)\vert=\vert d(h(x_{\delta}),h(x_0))-d(h(x_0),h(x_0))\vert=d(h(x_{\delta}),h(x_0))$ for every $\delta>\delta_0$.
That is, the net $\lbrace x_{\delta}\rbrace_{\delta\in \Delta}$ converges to $x_0$ in $t_p(\overline{G}^{M^X})$, which completes the proof.
\end{proof}

It is well known that the metric $\bar{d}:M\times M\rightarrow \R$ defined by $\bar{d}(m_1,m_2)\defi\min\lbrace d(m_1,m_2),1\rbrace$
 for all $m_1,m_2\in M$ induces the same topology as $d$. So, without loss of generality, we work with this metric from here on.

The following lemma reduces many questions related to
 a general metric space $M$ to the interval $[-1,1]$ (cf. \cite{Christensen1981}).

\blem\label{lemma_Christensen}
Let $X$ and $(M,d)$ be a topological and a metric space, respectively. If $G$ is a subset of $C(X,M)$,
then $G$ is equicontinuous at a point $x_0\in X$ if and only if $\nu(G\times K)$ is equicontinuous at it.
\elem
\begin{proof}
Assume that $G$ is equicontinuous at $x_0$. Given $\epsilon>0$, there is an open neighbouhood $U$ of $x_0$ such that
$d(g(x_0),g(x))<\epsilon$ for all $x\in U$ and $g\in G$. Let $\alpha\in K$, $x\in U$ and $g\in G$, then we have
$$\vert \nu(g,\alpha)(x_0)-\nu(g,\alpha)(x)\vert=\vert \alpha(g(x_0))-\alpha(g(x))\vert\leq d(g(x_0),g(x)) < \epsilon.$$

\noindent Conversely, assume that $\nu(G\times K)$ is equicontinuous in $x_0$. Given $\epsilon>0$, there is an open neighbouhood $U$ of $x_0$ such that $\vert \nu(g,\alpha)(x_0)-\nu(g,\alpha)(x)\vert<\epsilon$ for all $x\in U$, $g\in G$ and $\alpha\in K$.

For $g\in G$, consider the map $\alpha_{g(x_0)}\in K$. In order to finish the proof, it will suffice to observe that $$\vert \alpha_{g(x_0)}(g(x_0))-\alpha_{g(x_0)}(g(x))\vert=d(g(x),g(x_0))$$
\noindent for all $x\in U$ and $g\in G$.  
\end{proof}

\bcor\label{cor_Christensen}
Let $X$ and $(M,d)$ be a topological and a metric space, respectively, and let $G$ be an arbitrary subset of $C(X,M)$.
Then $G$ is (hereditarily) almost equicontinuous if and only if $\nu(G\times K)$ is (hereditarily) almost equicontinuous.

\ecor

\section{Proof of main results}

The following technical lemma is essential in most results along this paper.
The construction of the proof is based on an idea that
appears in \cite{Talagrand1979} and \cite{Cascales2000}.
We recall that a topological space is
\emph{hemicompact} if it has a sequence of compact subsets such that every compact
subset of the space lies inside some compact set in the sequence. Every compact space or
every locally compact and Lindel\"{o}f space is hemicompact.

\blem\label{Construccion_CANTOR}
Let $X$ and $(M,d)$ be a \v{C}ech-complete space and a hemicompact metric space, respectively,
and let $G$ be a subset of $C(X,M)$ such that $\overline{G}^{M^X}$ is compact.
If $G$ is not almost equicontinuous, then for every $G_{\delta}$ and dense subset $F$ of $X$ there exists
a countable subset $L$ in $G$, a compact separable subset $C_F\subseteq F$, a compact subset $N\subseteq M$
and a continuous and surjective map $\Psi$ of $C_F$ onto the Cantor set~$2^{\omega}$ such that
{for every $l\in L$ there exists a continuous map $l^*:2^{\omega}\rightarrow N$ satisfying that the following diagram is commutative}

\noindent \emph{Diagram 1:}
\[
\xymatrix{ C_F \ar@{>}[rr]^\Psi \ar[dr]^{l\vert_{C_F}} & & 2^\omega\ar[dl]_{l^*} \\ & N &}
\]

\noindent Furthermore, the subset $L^*\defi \{ l^* : l\in L\}\subseteq C(2^{\omega},N)$ separate points in $2^\omega$ and is not almost equicontinuous on $2^\omega$.

\elem
\begin{proof}
Let $F$ be a $G_{\delta}$ and dense subset of $X$. Then there is a sequence $\lbrace W_n\rbrace_{n=1}^{\infty}$ of open dense subsets of $X$ such that
$W_{s}\subseteq W_{r}$ if $r<s$ and $F=\bigcap\limits_{n=1}^{\infty}W_n$.

Let $\lbrace M_n\rbrace_{n<\omega}$ be a sequence of compact subsets, that we obtain by hemicompactness such that $M=\bigcup\limits_{n<\omega}M_n$ and for every compact subset $K\subseteq M$ there is $n<\omega$ such that $K\subseteq M_n$.

For each $n<\omega$ we consider the closed subset $X_n=\lbrace x\in X:g(x)\in M_n\quad \forall g\in G\rbrace$. We claim that $X=\bigcup\limits_{n<\omega}X_n$. Indeed, let $x\in X$. Since $\overline{G}^{M^X}\subseteq M^X$ is compact and the $x$th projection $\pi_x$ is continuous,
then $\pi_x(\overline{G}^{M^X})\subseteq M$ is compact. So, there is $n_x<\omega$ such that $\pi_x(\overline{G}^{M^X})\subseteq M_{n_x}$ by hemicompactness.
Therefore $x\in X_{n_x}$.

Since $G$ is not almost equicontinuous there exists a nonempty open subset $U$ of $X$ and $\epsilon>0$ such that for all nonempty  open subset $V\subseteq U$ there exists a function $g_V\in G$ such that $diam(g_V(V))\geq 2\epsilon>\gep$ by Lemma \ref{lem1}.

Note that $U$ is \v{C}ech-complete. If we express $U=\bigcup\limits_{n\in \omega}(U\cap X_n)$, by Baire's theorem, there is $n_0<\omega$ such that $\tilde{U}\defi int_U(U\cap X_{n_0})\neq \emptyset$ and open in $X$.

Set $C=\overline{\tilde{U}}^{X_{n_0}}$, which is closed in $X$, and $O_n=W_n\cap \tilde{U}=W_n\cap \tilde{U}\cap C$ that is open and dense in $C$ for each $n<\omega$. Then $O_{s}\subseteq O_{r}$ if $r<s$ and $H=\bigcap\limits_{n=1}^{\infty}O_n\subseteq F$ is a dense $G_{\delta}$ subset of $C$, which is a Baire space.
Remark further that $g(x)\in M_{n_0}$ for all $x\in C$ and $g\in G$.
Since $M_{n_0}$ is compact, every function $f\in C(C,M_{n_0})$ can be extended to a continuous function $f^{\beta}\in C(\beta C,M_{n_0})$.
Set $G^\gb=\{g^\gb : g\in G\}\subseteq C(\gb C,M_{n_0})$.

The space $C$, being \v{C}ech-complete, is a dense $G_{\delta}$ subset of its Stone-\v{C}ech compactification $\gb C$. Therefore, since $H$
is a $G_\gd$ subset of $C$, it follows that $H$ also is a dense $G_{\delta}$ subset of $\beta C$.
Consider a sequence $\lbrace E_n\rbrace_{n=1}^{\infty}$ of open dense subsets of
$\beta C$ such that $E_{s}\subseteq E_{r}$ if $r<s$ and $H=\bigcap\limits_{n=1}^{\infty}E_n$.
We have that $H=\bigcap\limits_{n=1}^{\infty}(E_n\cap  O^\gb_n)$, where $O^\gb_n=\gb C\setminus\overline{(C\setminus O_n)}^{\gb C}$
is open in $\gb C$ and $O^{\gb}_n\cap C=O_n$.

By induction on $n=\vert t \vert$ with $t\in 2^{(\omega)}$, we construct a family $\lbrace U_t:t\in 2^{(\omega)}\rbrace$ of nonempty open subsets of
$\gb C$ and a family of countable functions $L\defi\lbrace g_t:t\in 2^{(\omega)}\rbrace\subseteq G$, satisfying the following conditions for all $t\in 2^{(\omega)}$:
\begin{enumerate}[(i)]
\item $U_{\emptyset}\subseteq \overline{U_{\emptyset}}^{\gb C}\subseteq O^\gb_0\defi \gb C\setminus\overline{(C\setminus \tilde{U})}^{\gb C}$
(remark that $O^\gb _0\cap C=\tilde U$);
\item $U_{ti}\subseteq \overline{U_{ti}}^{\gb C}\subseteq E_{\vert t\vert}\cap O^\gb_{\vert t\vert}\cap U_t$ for $i=0,1$ (where $E_0\defi \beta C$);
\item $U_{t0}\cap U_{t1}=\emptyset$;
\item $d( g_t(x), g_t(y)) >\epsilon$, $\forall x\in U_{t0}\cap C$ and $\forall y\in U_{t1}\cap C$;
\item whenever $s,t\in 2^{(\omega)}$ and $\vert s\vert<\vert t \vert$, $diam(g_s(U_{tj}\cap C))<\frac{1}{\vert t\vert}$ for $j=0,1$.
\end{enumerate}
Indeed, if $n=0$, by regularity we can find $U_{\emptyset}$ an open set in $\gb C$ such that
$U_{\emptyset}\subseteq \overline{U_{\emptyset}}^{\gb  C}\subseteq E_0\cap O^\gb_0$.
For $n\geq 0$, suppose $\lbrace U_t:\vert t\vert \leq n\rbrace$ and $\lbrace g_t:\vert t\vert <n\rbrace$ have been constructed satisfying $(i)-(v)$.
Fix a $t\in 2^{(\omega)}$ with $\vert t\vert =n$.
Since $U_t$ is open in $\gb C$, there is an open set $A_t$ in $X$ such that $U_t\cap C=A_t\cap C$.
Therefore  
$$U_t\cap C=(A_t\cap C)\cap O^\gb_0=A_t\cap (O^\gb_0\cap C)=A_t\cap \tilde U$$ is open in $X$ and included in $U$.

By assumption there exist $g_t\in G$ such that $diam(g_t(U_t\cap C))>\epsilon$.
Consequently, we can find $x_t,y_t\in V_t\cap C$ such that $d( g_t(x_t),g_t(y_t)) >\epsilon$. By continuity, we can select two open
disjoint neighbourhoods in $\gb C$, $S_{t0}$ and $S_{t1}$ of $x_t$ and $y_t$, respectively, satisfying conditions $(iii)$ and $(iv)$.

If $i\in\lbrace 0,1\rbrace$, observe that $U_t\cap S_{ti}\cap O^\gb_0$ is open in $\gb C$ and nonempty.
Since $E_{\vert t\vert}\cap O^\gb_{\vert t\vert}$ is dense in $\gb C$ then
$U_t\cap S_{ti}\cap E_{\vert t\vert}\cap O^\gb_{\vert t\vert}$ is a nonempty open subset of $\gb C$. By regularity there exists a
 nonempty open subset $U_{ti}$ of $\gb C$ such that
 $U_{ti}\subseteq \overline{U_{ti}}^{\gb C}\subseteq U_t\bigcap S_{ti}\bigcap E_{\vert t\vert}\cap O^\gb_{\vert t\vert}$.
Therefore, $U_{t0}$ and $U_{t1}$ satisfies conditions $(ii), (iii)$ and $(iv)$ and, by continuity, we can adjust the open sets to satisfy $(v)$.

Set $K\defi\bigcap\limits_{n=0}^{\infty}\bigcup\limits_{\vert t \vert=n}\overline{U_t}^{\gb C}$, which is closed in $\gb C$ and,
as a consequence, also compact. Remark that we can express
$K=\bigcup\limits_{\sigma\in 2^{\omega}}\bigcap\limits_{n=0}^{\infty}\overline{U_{\sigma\vert n}}^{\gb C}$. Therefore, for each $\sigma\in 2^{\omega}$,
we have $\bigcap\limits_{n=0}^{\infty}\overline{U_{\sigma\vert n}}^{\beta C}\neq \emptyset$ by the compactness of $\gb C$, which implies $K\neq \emptyset$.
Furthermore,{ since $K\subseteq \bigcap\limits_{n=0}^{\infty}(E_n\cap O^\gb_n)=H\subseteq F$, it follows that $K$ is contained in $F$.}

Let $\Psi:K\rightarrow 2^{\omega}$ be the canonical map defined such that
$\Psi^{-1}(\sigma)=\bigcap\limits_{n=0}^{\infty}\overline{U_{\sigma\vert n}}^{\beta C}$ for all $\sigma\in 2^\omega$.
Clearly $\Psi$ is onto and continuous. Observe that for each $t\in 2^{(\omega)}$ and $\sigma\in 2^{\omega}$,
$g_t(\Psi^{-1}(\sigma))$ is a singleton by $(iv)$. Therefore, $g_t$ lifts to a continuous function $g_t^*$ on $2^{\omega}$
such that $g_t(x)=g_t^*(\Psi(x))$ for all $x\in K$.

Take a countable subset $D$ of $K$ such that $\Psi(D)=2^{(\omega)}$ and makes $\Psi\vert_D$ injective. Set $C_F\defi\overline{D}^{K}$. Note that $2^{(\omega)}$  is a countable dense subset of  $2^{\omega}$.

We have that $\Psi_{|C_F}:C_F\rightarrow 2^{\omega}$ is an onto and continuous map. We consider the set $L^*\subseteq C(2^{\omega},M_{n_0})$ defined by $L^*=\lbrace l^*:l\in L\vert_{C_F}\rbrace$ that makes the diagram 1 commutative.
We claim that $L^*$ separates points in $2^\omega$ and, as a consequence, defines its topology. Indeed, let $\sigma,\sigma'\in 2^\omega$ be two arbitrary points such that $\sigma\neq \sigma'$. Since $\Psi$ is an onto map there exist $x,y\in C_F$ such that $\sigma=\Psi(x)$ and $\sigma'=\Psi(y)$. Therefore, $x\in \bigcap\limits_{n=0}^{\infty}\overline{U_{\sigma\vert n}}^{\beta C}$ and $y\in \bigcap\limits_{n=0}^{\infty}\overline{U_{\sigma'\vert n}}^{\beta C}$. Since $\sigma\neq \sigma'$, there is $n_0\in\omega$ such that $\sigma\vert n_0=\sigma'\vert n_0$ and $\sigma(n_0+1)\neq \sigma'(n_0+1)$. Taking $t=\sigma\vert n_0$, then by $(iv)$ we know that $d(g_t(x),g_t(y))> \epsilon$. So, $g_t^*(\sigma)\neq g_t^*(\sigma')$.


On the other hand, by the commutativity of Diagram 1, and taking into account
how $L$ and $L^*$ have been defined, it is easily seen that $L^*$ is not almost equicontinuous on $2^\omega$
using Example \ref{ex_cantor2}.
\end{proof}
%


Applying {Corollary D of} \cite{Cascales2000} by Cascales, Namioka and Vera and Facts \ref{fact_0}, \ref{fact_1}, \ref{fact_2} and \ref{fact_AE},
next result follows easily.


\bprp\label{pr_CNV}
Let $X$ be a compact space, $(M,d)$ be a compact metric space and let $G$ be a subset of $C(X,M)$.
If $(X,t_p(\overline{G}^{M^X}))$ is Lindel\"of, then $G$ is hereditarily almost equicontinuous.
\eprp

Using Lemma \ref{Construccion_CANTOR}, the constraints in Proposition \ref{pr_CNV} can be relaxed as the following result shows.

\bprp\label{lem_Lindelof}
Let $X$ be a \v{C}ech-complete space, $(M,d)$ be a compact metric space and let $G$ be a subset of $C(X,M)$.
If there exists a dense $G_{\delta}$ subset $F\subseteq X$ such that
$(F,t_p(\overline{G}^{M^X}))$ is Lindel\"of, then $G$ is almost equicontinuous.
\eprp
\bpf
Reasoning by contradiction, suppose that $G$ is not almost equicontinuous.  By Lemma \ref{Construccion_CANTOR}
there exists a compact separable subset $C_F$ of $F$, a continuous onto map $\Psi:C_F\rightarrow 2^{\omega}$,
and a countable subset  $L$ of $G$ such that the subset $L^*\subseteq C(2^{\omega},M)$ defined by
$l^*(\Psi(x))=l(x)$ for all $x\in C_F$
separate points in $2^\omega$ and is not almost equicontinuous.

Let $K_F$ be the closure of $C_F$ in $F$ with respect to the initial topology generated by the maps in $L$.
Using a compactness argument, it follows that if $p\in K_F$ then there is $x_p\in C_F$ such that
$l(p)=l(x_p)$ for all $l\in L$.
Indeed, let $p\in K_F$. Then there is a net $\lbrace x_{\delta}\rbrace_{\delta\in \Delta}\subseteq C_F$ that $t_p(L)$-converges to $p$. Since $C_F$ is compact there is a subnet $\lbrace x_{\gamma}\rbrace_{\gamma\in \Gamma}$ such that converges to $x_0\in C_F$. Given $l\in L$, we know that $\lim\limits_{\gamma\in \Gamma}l(x_{\gamma})=l(x_0)$ because $l$ is continuous. Therefore, $l(x_0)=\lim\limits_{\gamma\in \Gamma}l(x_{\gamma})=l(p)$. Consequently, we can extend $\Psi$ to a map
$\Phi: K_F\rightarrow 2^{\omega}$ by $\Phi(p)=\Psi(x_p)$ for all $p\in K_F$.

Let's see that $\Phi$ is well-defined. Let $p\in K_F$, suppose that there are $x_p,\tilde{x}_p\in C_F$ such that $x_p\neq \tilde{x}_p$ and $l(p)=l(x_p)=l(\tilde{x}_p)$ for all $l\in L$. Since the Diagram 1 commutes, we know that $l^*(\Psi(x_p))=l^*(\Psi(\tilde{x}_p))$ for all $l^*\in L^*$. So, $\Psi(x_p)=\Psi(\tilde{x}_p)$ because $L^*$ separates points in $2^\omega$.\\ Observe that the following diagram is commutative

\noindent \emph{Diagram 2:}
\[
\xymatrix{ K_F \ar@{>}[rr]^\Phi \ar[dr]^{l\vert_{K_F}} & & 2^\omega\ar[dl]_{l^*} \\ & M &}
\]

\noindent
Certainly, let $p\in K_F$, then there is $x_p\in C_F$ such that $\Phi(p)=\Psi(x_p)$. Given $l\in L$, we have that $l(p)=l(x_p)=l^*(\Psi(x_p))=l^*(\Phi(p))$.\\
We claim that $\Phi:(K_F,t_p(L))\rightarrow (2^{\omega},t_p(L^*))$ is also continuous. Indeed, let $\lbrace h_{\delta}\rbrace_{\delta\in \Delta}\subseteq K_F$ a net that $t_p(L)$-converges to $h_0\in K_F$. For each $\delta\in \Delta$ there is $x_{\delta}\in C_F$ such that $\Phi(h_{\delta})=\Psi(x_{\delta})$ and $l(h_{\delta})=l(x_{\delta})$ for all $l\in L$. Analogously, there is $x_0\in C_F$ such that $\Phi(h_0)=\Psi(x_0)$ and $l(h_0)=l(x_0)$ for all $l\in L$.\\
Since $C_F$ is compact there is a subnet $\lbrace x_{\gamma}\rbrace_{\gamma\in \Gamma}$ such that converges to $\tilde{x}\in C_F$. Given $l\in L$, we know that $\lim\limits_{\gamma\in \Gamma}l(x_{\gamma})=l(\tilde{x})$ because $l$ is continuous. On the other hand, we also have that $\lim\limits_{\gamma\in \Gamma}l(x_{\gamma})=\lim\limits_{\gamma\in \Gamma}l(h_{\gamma})=l(h_0)=l(x_0)$. Therefore, $l(\tilde{x})=l(x_0)$ for all $l\in L$. So, $\Psi(\tilde{x})=\Psi(x_0)$ because $L^*$ separates points in $2^\omega$. The continuity follows because $\lim\limits_{\gamma\in \Gamma} \Phi(h_{\gamma})=\lim\limits_{\gamma\in \Gamma} \Psi(x_{\gamma})=\Psi(\tilde{x})=\Psi(x_0)=\Phi(h_0)$.

Now, since $K_F$ is $t_p(L)$-closed in $F$, it follows that it is also $t_p(\overline{G}^{M^X})$-closed in $F$.

By our initial assumption, we have that $F$ is $t_p(\overline{G}^{M^X})$-Lindel\"of, which implies that
also $K_F$ is $t_p(\overline{G}^{M^X})$-Lindel\"of.

We claim that $(2^{\omega},t_p(\overline{L^*}^{M^{2^{\omega}}}))$ is also Lindel\"of.
Indeed, it is enough to prove that $\Phi$ is continuous on $K_F$ when it is equipped with the
$t_p(\overline{G}^{M^X})$-topology and $2^{\omega}$ is equipped with the $t_p(\overline{L^*}^{M^{2^{\omega}}})$-topology.

Take a map $k\in \overline{L^*}^{M^{2^{\omega}}}$ and let $\lbrace l^*_{\gamma}\rbrace_{\gamma\in\Gamma}\subseteq L^*$
be a net converging to $k$ pointwise on $2^\omega$. Since $\overline{G}^{M^X}$ is compact, we may assume wlog that
$\lbrace l_{\gamma}\rbrace_{\gamma\in\Gamma}\subseteq L$ $t_p(X)$-converges to $h\in \overline{G}^{M^X}$.
Therefore, for each $x\in K_F$ we have that
$k(\Phi(x))=\lim\limits_{\gamma\in \Gamma} l^*_{\gamma}(\Phi(x))=\lim\limits_{\gamma\in \Gamma} l_{\gamma}(x)=h(x)$.
That is $k\circ \Phi=h$. Since $h$ is continuous on $K_F$, the continuity of $\Phi$ follows.

By Proposition \ref{pr_CNV}, this implies that $L^*$ is a hereditarily almost equicontinuous family on $2^\omega$,
 which is a contradiction.
\epf

\bprp\label{cor_Lindelof}
Let $X$ be a \v{C}ech-complete space, $(M,d)$ be a metric space and let $G$ be a subset of $C(X,M)$ such that $\overline{G}^{M^X}$ is compact.
If there exists a dense $G_{\delta}$ subset $F\subseteq X$ such that
$(F,t_p(\overline{G}^{M^X}))$ is Lindel\"of, then $G$ is almost equicontinuous.
\eprp
\bpf
Let $K$ and $\nu$ defined as in Remark \ref{nota_K}. Since 
$\nu (\overline{G}^{M^X}\times K)$ is a compact subset of $[-1,1]^X$, it follows that $\overline{\nu(G\times K)}^{[-1,1]^X}=\nu (\overline{G}^{M^X}\times K)$.

By Lemma \ref{lemma_homeo_sugerencia} we know that $(F,t_p(\nu(\overline{G}^{M^X}\times K)))$ is Lindel\"of . Now, applying
Proposition \ref{lem_Lindelof} to the subset $\nu (G\times K)\subseteq C(X,[-1,1])$, it follows that $\nu (G\times K)$ is almost equicontinuous.
Therefore, $G$ is almost equicontinuous by Corollary \ref{cor_Christensen}.
\epf
The following lemma is known. We refer to \cite[Cor. 3.5]{Ferrer2012} for its proof.

\blem\label{lem_Lindelof_METRIZABLE}
Let $X$ be a Lindel\"of space, $(M,d)$ be a metric space. If $G$ is an equicontinuous subset of $C(X,M)$,
then $\overline{G}^{M^X}$ is metrizable.
\elem

We are now in position of proving Theorem \ref{resultado_simpli}.

\begin{proof}[\textbf{Proof of Theorem \ref{resultado_simpli}}]
$(b)\Rightarrow (c)$ Since $(\overline{G}^{M^X})\vert_F$ is compact metric, it follows by Lemma \ref{lem_AlmostEquicontinuity} that
there is a dense subset $E$ such that $G$ is equicontinuous at the points in $E$ with respect to $X$.
Since $E$ is dense in $F$, which is dense in $X$, it follows that $E$ is also be dense in $X$. Moreover, if $Y$ denotes
the $G_{\delta}$ subset of equicontinuity points of $G$ in $X$, since $E\subseteq Y$, it follows that $Y$, the set of equicontinuity points of $G$ is a dense $G_\delta$-set in $X$.
Set $K\defi (\overline{G}^{M^X})$.
The equicontinuity of $G$ at the points in $Y$
combined with the density of $E\subseteq F$ in $Y$, implies that the map
$\Theta : K\vert_F\longrightarrow K\vert_Y$ defined by $\Theta(f\vert_F)\defi f\vert_Y$ is a homeomorphism of $K\vert_F$ onto $K\vert_Y$.   

By our initial assumption we have that
$K\vert_F$ is compact and metrizable, which yields the metrizability of $K\vert_Y$.
Thus, the evaluation map
$Eval :Y\lra C_\infty(K\vert_Y,M)$ is a well defined and continuous map. We know that $C_\infty(K\vert_Y,M)$ is
a separable space by \cite[Cor. 4.2.18]{Engelking1989}. Therefore $(Eval(Y),t_\infty(K\vert_Y))$ and $(Y,t_\infty(K\vert_Y) )$ are Lindel\"of spaces.
As a consequence $(Y,t_p(K\vert_Y))$ must be also Lindel\"of and we are done.

$(c)\Rightarrow (a)$ This implication is Proposition \ref{cor_Lindelof}

$(a)\Rightarrow (b)$ Suppose that $X$ is \v{C}ech-complete and hereditarily Lindel\"of.
By Lemma \ref{lem1}, the subset, $F$, of equicontinuity points
of $G$ is a dense $G_\gd$-set in $X$, which is a Lindel\"of space by our initial assumption.
Since $G$ is equicontinuous on $F$, Lemma \ref{lem_Lindelof_METRIZABLE}
implies that $(\overline{G}^{M^{X}})\vert_F$ must be metrizable.
\end{proof}
\mkp

The following result can be found in \cite[Prop. 2.5 and Section 5]{Glasner2006} in the setting of compact metric spaces.
Notwithstanding this, the proof given there can be adapted easily for \v{C}ech-complete and hereditarily Lindel\"of spaces,
as it is formulated in the next proposition. A sketch of the proof is included here for completeness sake.

\bprp\label{GMU}
Let $X$ be a hereditarily Lindel\"of space, $(M,d)$ is a metric space and $G\subseteq C(X,M)$.
If $H\defi\overline{G}^{M^X}$ is compact and hereditarily almost equicontinuous, then $H$ is metrizable.
\eprp
\begin{proof}
The symbol $C_{\infty}(H,M)$ denote the space $C(H,M)$ equipped with the uniform convergence topology. Consider the map $eval:X\rightarrow C_{\infty}(H,M)$ defined by $eval(x)[h]\defi h(x)$ for all $x\in X$ and $h\in H$.

By Proposition \ref{HAE_FRAG} $X$ is fragmented by $\rho_{G,d}$. Thus, for each nonempty subset $A$ of $X$ and for each
 $\epsilon>0$ there exists a nonempty open subset $U$ of $X$ such that $U\cap A\neq \emptyset$ and $diam(h(U\cap A))\leq \epsilon$ for all $h\in H$. Thus, $d_{\infty}$-$diam(eval(U\cap A))\leq \epsilon$.

We claim that $eval(X)$ is separable. Indeed, pick $\epsilon>0$. Let $\mathcal{A}$ be the collection of all open subsets $O$ of $X$ such that $eval(O)$ can be covered by countably many sets of diameter less than $\epsilon$. Since $X$ is hereditarily Lindel\"of  there is a countable subfamily $\mathcal{B}$ of $\mathcal{A}$ such that $\bigcup\limits_{A\in \mathcal{A}}A=\bigcup\limits_{B\in \mathcal{B}}B$. Take $V\defi \bigcup\limits_{A\in \mathcal{A}}A$. Observe that $V$ is the largest element of $\mathcal{A}$. Let's see that $A\defi X\setminus V$ is empty. Assume that $A\neq \emptyset$. Then there is a nonempty set $U$ of $X$ such that $U\cap A\neq \emptyset$ and $d_{\infty}$-$diam(eval(U\cap A))\leq \epsilon$.
Since $eval(U\cup V)=eval(U\cap A)\cup eval(V)$ we know that $eval(U\cup V)$ can be covered by countably many sets of diameter less than $\epsilon$. So, $U\cup V\in \mathcal{A}$ and we arrive to a contradiction because $U\cap (X\setminus V)\neq \emptyset$. Since $X=V\in \mathcal{A}$ and $\epsilon$ was arbitrary $eval(X)$ is separable.

There is a dense and countable subset $D$ of $eval(X)$. We know that $D$ separates points of $H$ because $eval(X)$ also separates points. Let $\Delta D:H\rightarrow M^D$ be the diagonal product. Since $\Delta D$ is an embedding and $M^D$ is metrizable we conclude that $H$ is metrizable.
\end{proof}

Next result is due basically to Namioka \cite[Lemma 2.1]{Namioka1987}. It can also be found in \cite[Lemma 6.4.]{Glasner2004a}, where the reference to
Namioka is acknowledged. Again, we include a sketch of the proof here for completeness sake.

\blem\label{lem_simplificacion}
Let $X$, $Y$ and $(M,d)$ be two arbitrary compact spaces and a metric space, respectively, and let $G$ be a subset of $C(Y,M)$. Suppose that $p : X\lra Y$ is a continuous onto map. Then $G\circ p \defi \{g\circ p : g\in G\}\subseteq C(X,M)$ is hereditarily almost equicontinuous if and only if $G$ is also hereditarily almost equicontinuous.
\elem
\begin{proof}
In order to prove this result, we will apply Lemma \ref{lem1}.
Assume that $G\circ p$ is hereditarily almost equicontinuous. Let $A$ be a closed (and compact) subset of $Y$, $U$ be a nonempty relatively open set in A and $\epsilon>0$. By Zorn's Lemma, there exists a minimal compact subset $Z$ of $X$ such that $p(Z)=A$. Since $\tilde{U}\defi p^{-1}(U)\cap Z$ is a nonempty relatively open set in $Z$ and $(G\circ p)\vert_Z$ is almost equicontinuous there is a nonempty relatively open set $\tilde{V}\subseteq \tilde{U}$ in $Z$ such that $diam((g\circ p)(\tilde{V}))< \epsilon$ for all $g\in G$. Let $V\defi A\setminus p(Z\setminus \tilde{V})$, that is relatively open set in $A$. We claim that $V\neq\emptyset$. Indeed, assume that $V=\emptyset$. Then $A=p(Z\setminus \tilde{V})$ and this contradicts the minimality of $Z$. Since $V\subseteq p(\tilde{V})$ we have that $diam(g(V))<\epsilon$ for all $g\in G$.

Conversely, let $Z$ be a closed subset of $X$, $\tilde{U}$ be a nonempty relatively open set in $Z$ and $\epsilon>0$. Consider the closed subset $W_0\defi \overline{p(\tilde{U})}$ of $Y$. Since $G\vert_{W_0}$ is almost equicontinuous there is a nonempty relatively open set $V_0$ in $Y$ such that $V_0\cap W_0\neq \emptyset$ and $diam(g(V_0\cap W_0))<\epsilon$ for all $g\in G$. Take $\tilde{V}\defi p^{-1}(V_0)\cap \tilde{U}$. Since $\tilde{V}$ is a nonempty relatively open set in $Z$ and $p(\tilde{V})\subseteq V_0\cap W_0$ we conclude that $diam(g(p(\tilde{V})))<\epsilon$ for all $g\in G$.
\end{proof}

\brem
If the map $p$ of the previous lemma is open or quasi-open we obtain the same result for almost equicontinuity. Recall that a map $f:X\rightarrow Y$ between two topological spaces is \textit{quasi-open} if for any nonempty open set $U\subseteq X$ the interior of $f(U)$ in $Y$ is nonempty.
\erem
\begin{proof}
Let $U$ be a nonempty open set of $Y$ and $\epsilon>0$. Since $G\circ p$ is almost equicontinuous and $\tilde{U}=p^{-1}(U)$ is an open subset of $X$ there is a nonempty open subset $\tilde{V}\subseteq \tilde{U}$ of $X$ such that $diam((g\circ p)(\tilde{V}))<\epsilon$ for all $g\in G$. Since the nonempty open set $V\defi int (p(\tilde{V}))$ is included in $p(\tilde{V})$ we have that $diam(g(V))<\epsilon$ for all $g\in G$.

Conversely, let $\tilde{U}$ be a nonempty open set of $X$ and $\epsilon>0$. Take $U\defi int(p(\tilde{U}))\neq\emptyset$. Since $G$ is almost equicontinuous there is a nonempty open subset $V\subseteq U$ of $Y$ such that $diam(g(V))<\epsilon$ for all $g\in G$. So, taking the open subset $\tilde{V}\defi p^{-1}(V)\cap \tilde{U}$, we conclude that $diam((g\circ p)(\tilde{V}))<\epsilon$ for all $g\in G$.
\end{proof}

%


\bprp\label{prop_a'_b'}
Let $X$ be a \v{C}ech-complete space, $(M,d)$ be a hemicompact metric space and $G\subseteq C(X,M)$ such that $\overline{G}^{M^X}$ is compact.
Then the following conditions are equivalent:
\begin{enumerate}[(a)]
\item $G$ is hereditarily almost equicontinuous.
\item $L$ is {hereditarily} almost equicontinuous on $F$, for all $L\in [G]^{\leq \omega}$ and $F$ a separable and compact subset of $X$.
\end{enumerate}
\eprp
\begin{proof}
$(a)$ implies $(b)$ is trivial. To see the other implication, assume, reasoning by contradiction, that (a) does not hold.
Then there must be some closed subset $A\subseteq X$
such that $G\vert_{A}$ is not almost equicontinuous.
By Lemma \ref{Construccion_CANTOR} there exists a compact and separable subset $F$ of $X$, an onto and continuous map $\Psi:F\rightarrow 2^{\omega}$, and a countable subset $L$ of $G$ such that the subset $L^*\subseteq C(2^{\omega},M)$ defined by $l^*(\Psi(x))=l(x)$ for all $x\in F$ is not almost equicontinuous. Therefore, $L$ is not hereditarily almost equicontinuous on $F$ by Lemma \ref{lem_simplificacion} and we arrive to a contradiction.
\end{proof}

We can now prove Theorem B.

\begin{proof}[\textbf{Proof of Theorem \ref{resultado_6}}]
$(b)\Rightarrow (a)$ is a direct consequence of Proposition \ref{prop_a'_b'} and Corollary \ref{cor_Christensen}.

$(a)\Rightarrow (c)$ Let $L\in[G]^{\leq \omega}$ and let $F$ be a separable and compact subset of $X$. $L$  defines an equivalence relation
on $F$  by $x\sim y$ if and only if $l(x)=l(y)$ for all $l\in L$. If $\widetilde{F}=F/{\sim}$ is the compact quotient space and
$p:F\rightarrow \widetilde{F}$ denotes the canonical quotient map,  each $l\in L$ has associated a map
$\tilde{l}\in C(\widetilde{F},M)$ defined as $\tilde{l}(\tilde{x})\defi l(x)$ for any $x\in F$ with $p(x)=\tilde{x}$.
Furthermore, if $\tilde{L}\defi \{\widetilde{l} : l\in L\}$, we can extend this definition to the closure of $\tilde{L}$ in $M^{\widetilde{F}}$.
Thus, each $l\in\overline{L}^{M^F}$ has associated a map $\tilde{l}\in\overline{\tilde{L}}^{M^{\widetilde{F}}}$ such that $\tilde{l}\circ p=l$.
By construction, we have that $\tilde{L}$ separates the points in $\widetilde{F}$. Since $\tilde{L}$ is countable it follows that $(\widetilde{F},t_p(\tilde{L}))$ is a compact metric space. On the other hand, $G$ is hereditarily almost equicontinuous
on $X$. Applying Lemma \ref{lem_simplificacion} to $F$ and $\widetilde{F}$, it follows that $\widetilde{L}$ is hereditarily almost equicontinuous on $\widetilde{F}$.
Therefore, the space $\overline{\tilde{L}}^{M^{\widetilde{F}}}$ is metrizable by Proposition \ref{GMU}. In order to finish the proof, it suffices
to remark that $\overline{L}^{M^F}$ is canonically homeomorphic to $\overline{\tilde{L}}^{M^{\widetilde{F}}}$ (see Fact \ref{lem_hom}).

$(c)\Rightarrow (d)$ Let $L\in[G]^{\leq \omega}$ and let $F$ be a separable and compact subset of $X$. We know that $H\defi ((\overline{L}^{M^X})\vert_F,t_p(F))$ is compact metric.
Since $F$ is separable, we have that $l(F)$ is a separable
for every $l\in L$. Hence $N\defi\overline{\bigcup\limits_{l\in L} l(F)}^M$ is a separable subset of $M$. Now, remark that $M$ can be replaced by $N$ without loss of generality. On the other hand, since $F\subseteq C(H,M)$ and $H$ is compact metric, it follows
that $(F,t_{\infty}(H))$ is separable and metrizable by \cite[Cor. 4.2.18]{Engelking1989}, which implies that $(F,t_{\infty}(H))$ is Lindel\"of.
Since the the topology $t_{p}(H)$ is weaker than $t_{\infty}(H)$, we deduce that $(F,t_{p}(H))$ must be Lindel\"{o}f.

$(d) \Rightarrow (b)$ By Lemma \ref{lemma_homeo_sugerencia}, for all $L\in [G]^{\leq \omega}$
and $F$ a separable compact subset of $X$, we have that $(F,t_p(\nu(\overline{L}^{M^X}\times K)))$ is Lindel\"of . Applying  \cite[Corollary D]{Cascales2000}, it follows that $\nu(\overline{L}^{M^X}\times K)$ is hereditarily almost equicontinuous for all $L\in [G]^{\leq \omega}$ and $F$ a separable compact subset of $X$.
Thus, Corollary \ref{cor_Christensen} yields $(b)$.
\end{proof}

\section{Appendix}
It is well known that for every compact metric space $(M,d)$, there is a canonical continuous one-to-one mapping
$\mathcal{E}_M : M\lra [0,1]^{\omega}$ that embeds $M$ into $[0,1]^{\omega}$ as a closed subspace.
Let $\rho_n : [-1,1]\lra [0,1]$ the map defined by $\rho_n(r)= \frac{|r|}{2^n}$ for every $n<\omega$.
Along this paper, we will consider that $[0,1]^{\omega}$ is equipped with the metric $\rho$ defined by
\[\rho((x_n),(y_n))=\sum_{n<\omega} \rho_n(x_n-y_n)\]
\bkp
The proof of the following lemma is obtained by a standard argument of compactness, using the continuity of
$\mathcal{E}_M^{-1}$ and that every continuous map defined on a compact space is uniformly continuous.
We omit its proof here.

\bfct\label{fact_0}
Let $(M,d)$ be a compact metric space. Let $\mathcal{E}_M : M\lra [0,1]^{\omega}$ denote
its attached embedding into $[0,1]^{\omega}$, and let $\pi_n:[0,1]^{\omega}\rightarrow [0,1]$ denote the \emph{$n$th} canonical projection.
Then, for every $\epsilon >0$, there is $\delta >0$ and $n_0<\omega$ such that if $(x,y)\in M\times M$ and
$\rho_n(\pi_n(\mathcal{E}_M(x))- \pi_n(\mathcal{E}_M(y)))<\delta/2n_0$ for $n\leq n_0$ then $d(x,y)<\epsilon$.
\efct

We know recall some simple remarks that will be used along the paper. 

\bfct\label{fact_1}
Let $X$ be a topological space and $(M,d)$ a compact metric space. If $\pi_n$ is the $n$th projection mapping defined above,
then the following map is continuous if we consider that the two spaces have the topology of pointwise convergence.
\[\pi_n^*:M^X \rightarrow[0,1]^X\]
defined by $\pi_n^*(f)\defi  \pi_n\circ \mathcal{E}_M\circ f$, $f\in M^X$, for each $n<\omega$.
\efct

For each $S\subseteq M^X$ and each $n<\omega$ we define $S_n\defi\pi_n^*(S)$.

\bfct\label{fact_2}
{Let $X$ be a Baire space, $(M,d)$ be a compact metric space, $G\subseteq C(X,M)$ and $H\defi\overline{G}^{M^X}$.} Then $H_n=\overline{G_n}^{[0,1]^X}$.
\efct
\bpf
Indeed, since $\pi_n^*$ is continuous we have that
$H_n=\pi_n^*(H)=\pi_n^*(\overline{G}^{M^X})\subseteq \overline{\pi_n^*(G)}^{[0,1]^X}=\overline{G_n}^{[0,1]^X}$.
For the reverse inclusion, remark that $\overline{G_n}^{[0,1]^X}$ is the {smallest} closed subset that contains $G_n$ and $G_n\subseteq H_n$.
\epf


\bfct\label{fact_AE}
Let $X$ be a Baire space, $(M,d)$ be a compact metric space and $G\subseteq C(X,M)$. If $G_n$ is almost equicontinuous for every $n<\omega$, then $G$ is almost equicontinuous.
\efct
\begin{proof}
For each $n\in \omega$ there exists a dense $G_{\delta}$ subset $D_n$ of $X$ such that $G_n$ is equicontinuous on $D_n$. Since $X$ is a Baire space, the $D=\bigcap\limits_{n<\omega}D_n$ is dense in $X$. We claim that $G$ is equicontinuous in $D$. Indeed, let $x_0\in D$ and $\epsilon>0$. By Fact \ref{fact_0} we get $\delta>0$ and $n_0<\omega$. Take $\epsilon_0=\frac{\delta}{2n_0}$. For each $n<n_0$, being $G_n$ equicontinuous in $x_0$, there is an open neighbourhood $U_n$ of $x_0$ such that $\vert g_n(x_0)-g_n(x)\vert<\epsilon_0$ for all $x\in U_n$ and $g_n\in G_n$. Consider the open neighbourhood $U=\bigcap\limits_{n<n_0}U_n$ of $x_0$.
So, let an arbitrary $g\in G$ and $x\in U$, then $\rho_n(\pi_n(\mathcal{E}_M(g(x_0)))- \pi_n(\mathcal{E}_M(g(x))))=\rho_n(\pi^*_n(g)(x_0)-\pi^*_n(g)(x))=\frac{\vert \pi^*_n(g)(x_0)-\pi^*_n(g)(x) \vert }{2^n}<\frac{\epsilon_0}{2^n}\leq\frac{\delta}{2n_0}$. Consequently, $d(g(x_0),g(x))<\epsilon$ by Fact \ref{fact_0}.
\end{proof}


\bfct\label{fact_3}
The diagonal map $\Delta:H\rightarrow \prod\limits_{n<\omega} H_n$ defined by $\Delta(h)=(\pi_n\circ \mathcal{E}_M\circ h)_{n<\omega}$ for each $h\in H$, is a homeomorphism of $H$ onto its image.
\efct
\bfct
 Given a subset $L\subseteq G$, it defines an equivalence relation on $X$  by $x\sim y$ if and only if $l(x)=l(y)$ for all $l\in L$.
 Let $\widetilde{X}=X/{\sim}$ be the quotient space and let $p:X\rightarrow \widetilde{X}$ denote the canonical quotient map,
then each $l\in L$ has associated a map
$\tilde{l}\in C(\widetilde{X},M)$ defined as $\tilde{l}(\tilde{x})\defi l(x)$ for any $x\in X$ with $p(x)=\tilde{x}$.
Furthermore, if $\tilde{L}\defi \{\widetilde{l} : l\in L\}$, we can extend this definition to the closure of $\tilde{L}$.
Thus, each $l\in\overline{L}^{M^X}$ has associated a map $\tilde{l}\in\overline{\tilde{L}}^{M^{\widetilde{X}}}$ such that $\tilde{l}\circ p=l$.
\efct


\bfct\label{lem_hom}
Let $L$ be a countable subset of $G\subseteq C(X,M)$. We denote by $X_L$ the topological space $(\widetilde{X},t_p(\tilde{L}))$, which is metrizable because $\tilde{L}$ is countable.
Consider the  map $
p^*:(M^{\tilde{X}},t_p(\tilde{X})) \rightarrow (M^X,t_p(X))$ defined by $p^*(\tilde{f})=\tilde{f}\circ p$, for each $\tilde{f}\in M^{\tilde{X}}$.
Then $p^*$ is a homeomorphism of $\overline{\tilde{L}}^{M^{\tilde{X}}}$ onto $\overline{L}^{M^X}$.
\efct
\begin{proof}
 We observe that $p^*$ is continuous, since a net $\lbrace \tilde{f_{\alpha}}\rbrace_{\alpha\in A}$ $t_p(\tilde{X})$-converges to $\tilde{f}$ in $\overline{\tilde{L}}^{M^{\tilde{X}}}$ if and only if $\lbrace \tilde{f_{\alpha}}\circ p\rbrace_{\alpha\in A}$ $t_p(X)$-converges to $\tilde{f}\circ p$ in $\overline{L}^{M^{X}}$.\\
  Let's see that $p^*(\overline{\tilde{L}}^{M^{\tilde{X}}})=\overline{L}^{M^X}$. Indeed, since $p^*$ is continuous we have that $p^*(\overline{\tilde{L}}^{M^{\tilde{X}}})\subseteq \overline{p^*(\tilde{L})}^{M^X}=\overline{L}^{M^X}$. We have the other inclusion because $\overline{L}^{M^X}$ is the smaller closed set that contains $L$ and $L\subseteq p^*(\overline{\tilde{L}}^{M^{\tilde{X}}})$.\\
 Let $\tilde{f},\tilde{g}\in \overline{\tilde{L}}^{M^{\tilde{X}}}$ such that $\tilde{f}\neq\tilde{g}$. Then there exists $\tilde{x}\in \tilde{X}$ such that $\tilde{f}(\tilde{x})\neq\tilde{g}(\tilde{x})$. Let $x\in X$ an element such that $\tilde{x}=p(x)$. Thus $(\tilde{f}\circ p)(x)\neq (\tilde{g}\circ p)(x)$. So, $p^*$ is injective because $\tilde{f}\circ p\neq \tilde{g}\circ p$.\\
Finally, we arrive to the conclusion that $p^*\vert_{\overline{\tilde{L}}^{M^{\tilde{X}}}}$ is a homeomorphism because it is defined between compact spaces.
\end{proof}

\section{Acknowledgments}
We are very grateful to the referee for a thorough report that helped considerably in improving the presentation of this paper.


\bibliography{BiblioResultados}
\bibliographystyle{plain}

\bigskip

\end{document}